\documentclass[a4paper,twosided, 11.7pt]{amsart}

\usepackage{amsfonts}
\usepackage{amssymb}
\usepackage{latexsym}
\usepackage{graphics}
\usepackage[2emode]{psfrag}
\usepackage{mathrsfs}
\usepackage{amsthm}
\usepackage{amsmath}
\usepackage[all]{xy}
\usepackage{enumerate}
\usepackage[latin1]{inputenc}
\usepackage{lscape}
\usepackage{a4wide}

\xyoption{2cell} \xyoption{curve} \xyoption{rotate} \xyoption{v2}

\begin{document}

\def\id{\textrm{{\small 1}\normalsize\!\!1}}
\def\To{{\multimap\!\to}}
\def\bigperp{{\LARGE\textrm{$\perp$}}}
\newcommand{\QED}{\hspace{\stretch{1}}
\makebox[0mm][r]{$\Box$}\\}

\def\AA{{\mathbb A}}
\def\BB{{\mathbb B}}
\def\CC{{\mathbb C}}
\def\DD{{\mathbb D}}
\def\EE{{\mathbb E}}
\def\FF{{\mathbb F}}
\def\GG{{\mathbb G}}
\def\HH{{\mathbb H}}
\def\II{{\mathbb I}}
\def\JJ{{\mathbb J}}
\def\KK{{\mathbb K}}
\def\LL{{\mathbb L}}
\def\MM{{\mathbb M}}
\def\NN{{\mathbb N}}
\def\OO{{\mathbb O}}
\def\PP{{\mathbb P}}
\def\QQ{{\mathbb Q}}
\def\RR{{\mathbb R}}
\def\SS{{\mathbb S}}
\def\TT{{\mathbb T}}
\def\UU{{\mathbb U}}
\def\VV{{\mathbb V}}
\def\WW{{\mathbb W}}
\def\XX{{\mathbb X}}
\def\YY{{\mathbb Y}}
\def\ZZ{{\mathbb Z}}

\def\aa{{\mathfrak A}}
\def\bb{{\mathfrak B}}
\def\cc{{\mathfrak C}}
\def\dd{{\mathfrak D}}
\def\ee{{\mathfrak E}}
\def\ff{{\mathfrak F}}
\def\gg{{\mathfrak G}}
\def\hh{{\mathfrak H}}
\def\ii{{\mathfrak I}}
\def\jj{{\mathfrak J}}
\def\kk{{\mathfrak K}}
\def\ll{{\mathfrak L}}
\def\mm{{\mathfrak M}}
\def\nn{{\mathfrak N}}
\def\oo{{\mathfrak O}}
\def\pp{{\mathfrak P}}
\def\qq{{\mathfrak Q}}
\def\rr{{\mathfrak R}}
\def\ss{{\mathfrak S}}
\def\tt{{\mathfrak T}}
\def\uu{{\mathfrak U}}
\def\vv{{\mathfrak V}}
\def\ww{{\mathfrak W}}
\def\xx{{\mathfrak X}}
\def\yy{{\mathfrak Y}}
\def\zz{{\mathfrak Z}}

\def\aaa{{\mathfrak a}}
\def\bbb{{\mathfrak b}}
\def\ccc{{\mathfrak c}}
\def\ddd{{\mathfrak d}}
\def\eee{{\mathfrak e}}
\def\fff{{\mathfrak f}}
\def\ggg{{\mathfrak g}}
\def\hhh{{\mathfrak h}}
\def\iii{{\mathfrak i}}
\def\jjj{{\mathfrak j}}
\def\kkk{{\mathfrak k}}
\def\lll{{\mathfrak l}}
\def\mmm{{\mathfrak m}}
\def\nnn{{\mathfrak n}}
\def\ooo{{\mathfrak o}}
\def\ppp{{\mathfrak p}}
\def\qqq{{\mathfrak q}}
\def\rrr{{\mathfrak r}}
\def\sss{{\mathfrak s}}
\def\ttt{{\mathfrak t}}
\def\uuu{{\mathfrak u}}
\def\vvv{{\mathfrak v}}
\def\www{{\mathfrak w}}
\def\xxx{{\mathfrak x}}
\def\yyy{{\mathfrak y}}
\def\zzz{{\mathfrak z}}

\newcommand{\aA}{\mathscr{A}}
\newcommand{\bB}{\mathscr{B}}
\newcommand{\cC}{\mathscr{C}}
\newcommand{\dD}{\mathscr{D}}
\newcommand{\eE}{\mathscr{E}}
\newcommand{\fF}{\mathscr{F}}
\newcommand{\gG}{\mathscr{G}}
\newcommand{\hH}{\mathscr{H}}
\newcommand{\iI}{\mathscr{I}}
\newcommand{\jJ}{\mathscr{J}}
\newcommand{\kK}{\mathscr{K}}
\newcommand{\lL}{\mathscr{L}}
\newcommand{\mM}{\mathscr{M}}
\newcommand{\nN}{\mathscr{N}}
\newcommand{\oO}{\mathscr{O}}
\newcommand{\pP}{\mathscr{P}}
\newcommand{\qQ}{\mathscr{Q}}
\newcommand{\rR}{\mathscr{R}}
\newcommand{\sS}{\mathscr{S}}
\newcommand{\tT}{\mathscr{T}}
\newcommand{\uU}{\mathscr{U}}
\newcommand{\vV}{\mathscr{V}}
\newcommand{\wW}{\mathscr{W}}
\newcommand{\xX}{\mathscr{X}}
\newcommand{\yY}{\mathscr{Y}}
\newcommand{\zZ}{\mathscr{Z}}

\newcommand{\Aa}{\mathcal{A}}
\newcommand{\Bb}{\mathcal{B}}
\newcommand{\Cc}{\mathcal{C}}
\newcommand{\Dd}{\mathcal{D}}
\newcommand{\Ee}{\mathcal{E}}
\newcommand{\Ff}{\mathcal{F}}
\newcommand{\Gg}{\mathcal{G}}
\newcommand{\Hh}{\mathcal{H}}
\newcommand{\Ii}{\mathcal{I}}
\newcommand{\Jj}{\mathcal{J}}
\newcommand{\Kk}{\mathcal{K}}
\newcommand{\Ll}{\mathcal{L}}
\newcommand{\Mm}{\mathcal{M}}
\newcommand{\Nn}{\mathcal{N}}
\newcommand{\Oo}{\mathcal{O}}
\newcommand{\Pp}{\mathcal{P}}
\newcommand{\Qq}{\mathcal{Q}}
\newcommand{\Rr}{\mathcal{R}}
\newcommand{\Ss}{\mathcal{S}}
\newcommand{\Tt}{\mathcal{T}}
\newcommand{\Uu}{\mathcal{U}}
\newcommand{\Vv}{\mathcal{V}}
\newcommand{\Ww}{\mathcal{W}}
\newcommand{\Xx}{\mathcal{X}}
\newcommand{\Yy}{\mathcal{Y}}
\newcommand{\Zz}{\mathcal{Z}}

\numberwithin{equation}{section}
\renewcommand{\theequation}{\thesection.\arabic{equation}}
\newcommand{\bara}[1]{\overline{#1}}
\newcommand{\ContFunt}[2]{\bara{\mathrm{Funt}}(#1,\,#2)}
\newcommand{\lBicomod}[2]{{}_{#1}\mM^{#2}}
\newcommand{\rBicomod}[2]{{}^{#1}\mM_{#2}}
\newcommand{\Bicomod}[2]{{}^{#1}\mM^{#2}}
\newcommand{\lrBicomod}[2]{{}^{#1}\mM^{#2}}
\newcommand{\rcomod}[1]{\mM^{#1}}
\newcommand{\rmod}[1]{\mM_{#1}}
\newcommand{\bimod}[1]{{}_{#1}\mM_{#1}}
\newcommand{\Bimod}[2]{{}_{#1}\mM_{#2}}
\newcommand{\Sf}[1]{\mathsf{#1}}
\newcommand{\Sof}[1]{S^{{\bf #1}}}
\newcommand{\Tof}[1]{T^{{\bf #1}}}
\renewcommand{\hom}[3]{\mathrm{Hom}_{#1}\left(\underset{}{}#2,\,#3\right)}
\newcommand{\Coint}[2]{\mathrm{Coint}(#1,#2)}
\newcommand{\IntCoint}[2]{\mathrm{InCoint}(#1,#2)}
\newcommand{\Coder}[2]{\mathrm{Coder}(#1,#2)}
\newcommand{\IntCoder}[2]{\mathrm{InCoder}(#1,#2)}
\newcommand{\lr}[1]{\left(\underset{}{} #1 \right)}
\newcommand{\Ext}[3]{\mathrm{Ext}_{\eE}^{#1}\lr{#2,\,#3}}
\newcommand{\equalizerk}[2]{\mathfrak{eq}_{#1,\,#2}^k}
\newcommand{\equalizer}[2]{\mathfrak{eq}_{#1,\,#2}}
\newcommand{\coring}[1]{\mathfrak{#1}}
\newcommand{\tensor}[1]{\otimes_{#1}}
\newcommand{\Mono}[1]{\rR_{(\coring{#1}\,:\,A)}}
\newcommand{\Rtensor}[1]{\underset{(\coring{#1}\,:\,A)}{\otimes}}
\newcommand{\cotensor}[1]{\square_{#1}}
\newcommand{\tensorbajo}[1]{\underset{#1}{\otimes}}
\newcommand{\td}[1]{\widetilde{#1}}
\newcommand{\tha}[1]{\widehat{#1}}
\newcommand{\wrcomod}[2]{ \mM^{#1}_{#2}}
\newcommand{\bd}[1]{\boldsymbol{#1}}
\newcommand{\tensork}{\otimes_{\Bbbk}}


\def\units{{\mathbb G}_m}
\def\rightact{\hbox{$\leftharpoonup$}}
\def\leftact{\hbox{$\rightharpoonup$}}

\def\*C{{}^*\hspace*{-1pt}{\Cc}}

\def\text#1{{\rm {\rm #1}}}

\def\smashco{\mathrel>\joinrel\mathrel\triangleleft}
\def\cosmash{\mathrel\triangleright\joinrel\mathrel<}

\def\ol{\overline}
\def\ul{\underline}
\def\dul#1{\underline{\underline{#1}}}
\def\Nat{\dul{\rm Nat}}
\def\Set{\dul{\rm Set}}

\renewcommand{\subjclassname}{\textup{2000} Mathematics Subject
     Classification}

\newtheorem{proposition}{Proposition}[section]
  \newtheorem{lemma}[proposition]{Lemma}
  \newtheorem{aussage}[proposition]{}
  \newtheorem{corollary}[proposition]{Corollary}
  \newtheorem{theorem}[proposition]{Theorem}

  \theoremstyle{definition}
  \newtheorem{definition}[proposition]{Definition}
  \newtheorem{example}[proposition]{Example}

  \theoremstyle{remark}
  \newtheorem{remark}[proposition]{Remark}

\title[Compatibility condition between ring and coring.]
{Compatibility condition between ring and coring.}
\date{\today}
\author{L. El Kaoutit}
\address{Departamento de \'Algebra. Facultad de Educaci\'on y Humanidades de Ceuta.
Universidad de Granada. El Greco Nº 10. E-51002 Ceuta, Spain}
\email{kaoutit@ugr.es}

\keywords{ (Co)Rings. Bi-(Co)modules. Monoidal categories. (Co)Monoides. (Co)Wreath.  Double Distributive law.\\
Research supported by grant MTM2004-01406 from the Ministerio de
Educaci\'{o}n y Ciencia of Spain. } 

\baselineskip 17pt

\begin{abstract}
We introduce the notion of bi-monoid in general monoidal category
generalizing by this the notion of bialgebra. In the case of
bimodules over a noncommutative algebra, we obtain a compatibility
condition between ring and coring whenever both structures admit the
same underlying bimodule.
\end{abstract}

\maketitle

\section*{Introduction}

Let $R$ be an associative algebra over a commutative ring with
identity $\Bbbk$, and consider its category of unital and
$\Bbbk$-central bimodule ${}_R\mM_R$ as a monoidal category with
multiplication the two variables functor defined by the tensor
product $-\tensor{R}-$ and with identity object the regular bimodule
${}_RR_R$. A compatibility condition between a ring structure
(monoid in $\bimod{R}$, see below) and coring structure (comonoid in
$\bimod{R}$, see below), is a well known problem in noncommutative
algebra. Even though, the object is the same (i.e., the underlying
$R$-bimodules structures coincide), there is no obvious way in which
a tensor product of bimodules can be equipped with a monoid
structure. In this direction M. E. Sweedler \cite[\S
5]{Sweedler:1975b} introduced the $\times_A$-product, and few years
later M. Takeuchi gave in \cite{Takeuchi:1977} the notion of
$\times_A$-bialgebra with noncommutative basering. This notion can
be seen as an approach to a compatibility condition problem by
taking, in appropriate way, a comonoid in $\bimod{R}$ and a monoid
in $\bimod{R\tensork R^o}$ ($R^o$ is the opposite ring of $R$).

In this note we give another approach to the compatibility condition
problem by considering a comonoid and monoid in the same monoidal
category. The basic ideas behind our approach are the notions of
wreath and cowreath recently introduced by S. Lack and R. Street in
\cite{Lack/Street:2002}, which are a generalization of distributive
law due to J. Beck \cite{Beck:1969}. In the bimodules category,
wreath and cowreath lead in a formal way to endow certain tensor
product within a structure of ring and of coring. A double
distributive law is then a wreath and cowreath induced,
respectively, by a ring and coring taking the same underlying
bimodule. In this way we arrive to the compatibility condition by
assuming the existence of a double distributive law for a bimodule
which admits a structures of both ring and coring.

In Section \ref{ELM} we review the Eilenberg-Moore categories
attached to a monoid (resp. comonoid) in a strict monoidal category.
We give, as in the case of bimodules category \cite{Kaoutit:2007pr},
a simplest and equivalent definition for a wreath (resp. cowreath)
over monoid (resp. comonoid), see Proposition \ref{A-Def-equ} (resp.
Proposition \ref{Def-equ}). In particular we show that the wreath
(resp. cowreath) product satisfies a universal property, Proposition
\ref{Univ-Prp} (resp. Proposition \ref{Univ-coPrp}). In Section
\ref{BiMon}, we use the notion of double distributive law
(Definition \ref{doble-Dl}) in order to prove an equivalent
compatibility conditions for an object with a structures of both
monoid and comonoid, Proposition \ref{Bimonoid}. An object
satisfying the equivalent conditions of such proposition is called a
bi-monoid. Section \ref{aplication} presents an application of
results stated in previous sections to the bimodules category. In
particular, we give an example which shows that the class of
bialgebras is "strictly contained" in the class of bi-monoid in the
monoidal category of $\Bbbk$-modules.

\bigskip
{\textsc{Notations and Basic Notions}:} Given any Hom-set category
$\cC$, the notation $X \in \aA$ means that $X$ is an object of
$\cC$. The identity morphism of $X$ will be denoted by $X$ itself.
The set of all morphisms $f:X \to X'$ in $\cC$, is denoted by
$\hom{\cC}{X}{X'}$. Let $\Mm$ be a strict monoidal category with
multiplication $-\tensor{}-$ and identity object $\II$. Recall from
\cite[\S VII]{Maclane:livro} that a comonoid in $\Mm$ is a
three-tuple $(\CC,\Delta,\varepsilon)$ consisting of an object $\CC$
and two morphisms $\Delta: \CC \to \CC\tensor{}\CC$
(comultiplication), $\varepsilon: \CC \to \II$ (counit) such that
$(\varepsilon\tensor{}\CC) \circ \Delta\,=\,\CC\,=\,
(\CC\tensor{}\varepsilon) \circ \Delta$ and $(\Delta\tensor{}\CC)
\circ \Delta\,=\, (\CC\tensor{}\Delta) \circ \Delta$. A morphisms of
comonoids $\phi: (\CC,\Delta,\varepsilon) \to
(\DD,\Delta',\varepsilon')$ is a morphism $\phi: \CC \to \DD$ in
$\Mm$ such that $\varepsilon' \circ \phi \,=\, \varepsilon$
(counitary property) and $\Delta' \circ \phi\,=\,
(\phi\tensor{}\phi) \circ \Delta$ (coassociativity property).
Dually, a monoid in $\Mm$ is a three-tuple $(\AA,\mu,\eta)$
consisting of an object $\AA \in \Mm$ and two morphisms $\mu:
\AA\tensor{}\AA \to \AA$ (multiplication), $\eta :\II \to \AA$
(unit) such that $\mu \circ (\eta\tensor{}\AA) \,=\, \AA\,=\, \mu
\circ (\AA\tensor{}\eta)$ and $ \mu \circ (\AA\tensor{}\mu)\,=\, \mu
\circ (\mu\tensor{}\AA)$. A morphism of monoid $\psi: (\AA,\mu,\eta)
\to (\TT,\mu',\eta')$ is a morphism $\psi: \AA \to \TT$ such that
$\psi \circ \eta\,=\, \eta'$ and $\psi \circ \mu\,=\, \mu' \circ
(\psi\tensor{}\psi)$.

A right $\CC$-comodule is a pair $(X,\rho^{X})$ consisting of an
object $X \in \Mm$ and a morphism $\rho^X: X \to X\tensor{}\CC$ (a
right $\CC$-coaction) such that $(X\tensor{}\varepsilon) \circ
\rho^X\,=\, X$ and $(\rho^{X}\tensor{}\CC) \circ \rho^X\,=\,
(X\tensor{}\Delta) \circ \rho^X$. A morphism of right
$\CC$-comodules $f:(X,\rho^X) \to (X',\rho^{X'})$ is a morphism
$f:X\to X'$ in the category $\Mm$ such that $ \rho^{X'} \circ f\,=\,
(f\tensor{}\CC) \circ \rho^X$. Left $\CC$-comodules and their
morphisms are similarly defined; we use the Greek letter
$\lambda^{-}$ to denote their left $\CC$-coactions. A
$\CC$-bicomodule is a three-tuple $(X,\rho^X,\lambda^X)$ where
$(X,\rho^X)$ is a right $\CC$-comodule and $(X,\lambda^X)$ is a left
$\CC$-comodule such that $ (\lambda^X\tensor{}\CC) \circ \rho^X\,=\,
(\CC\tensor{}\rho^X) \circ \lambda^X$. A morphism of
$\CC$-bicomodules is a morphism of left and of right
$\CC$-comodules. We use the notation $\hom{\CC-\CC}{-}{-}$ for the
sets of all $\CC$-bicomodules morphisms. Dually, a right
$\AA$-module is a pair $(P,\Sf{r}_P)$ consisting of an object $P \in
\Mm$ and a morphisms $\Sf{r}_P: P\tensor{}\AA \to P$ (a right
$\AA$-action) such that $\Sf{r}_P \circ (P\tensor{}\eta)\,=\, P$ and
$\Sf{r}_P \circ (P\tensor{}\mu)\,=\, \Sf{r}_P \circ
(\Sf{r}_P\tensor{}\AA)$. A morphism of right $\AA$-module $g:
(P,\Sf{r}_P) \to (P',\Sf{r}_{P'})$ is a morphism $g:P \to P'$ in
$\Mm$ such that $ g \circ \Sf{r}_P\,=\, \Sf{r}_{P'} \circ (g
\tensor{}\AA)$.
 Left $\AA$-module and their morphisms are similarly defined and we use the letter
$\Sf{l}_{-}$ to denote a left $\AA$-actions. An $\AA$-bimodule is a
three-tuple $(P,\Sf{r}_P,\Sf{l}_P)$ where $(P,\Sf{r}_P)$ is a right
$\AA$--module and $(P,\Sf{l}_P)$ is a left $\AA$--module such that $
\Sf{r}_P \circ (\Sf{l}_P\tensor{}\AA)\,=\, \Sf{l}_P \circ
(\AA\tensor{}\Sf{r}_P)$. A morphisms of $\AA$-bimodules is a
morphism of left and of right $\AA$-modules. We denote by
$\hom{\AA-\AA}{-}{-}$ the sets of all $\AA$-bimodules morphisms. For
a details on comodules corings, definitions and basic properties of
bicomodules over corings, the reader is referred to monograph
\cite{Brzezinski/Wisbauer:2003}.

\section{Review on (right) Eilenberg-Moore monoidal categories.}\label{ELM}

Let $\Mm$ denote a strict monoidal category with multiplication
$-\tensor{}-$ and identity object $\II$. In this section we review
the Eilenberg-Moore monoidal categories \cite{Lack/Street:2002}
associated to a monoid and to a comonoid both defined in $\Mm$. We
start by considering a comonoid $\CC$ in $\Mm$ with a structure
morphisms $\Delta: \CC \to \CC \tensor{}\CC$ and $\varepsilon: \CC
\to \II$.

\subsection{The monoidal category $\rR_{\CC}^{\,c}$}\label{RC} $ $ \\
This is the right Eilenberg-Moore monoidal category associated to
the comonoid $\CC$ (see \cite{Lack/Street:2002} for a general
context) defined as follows

\paragraph{\emph{Objects of}  $\rR_{\CC}^{\,c}$:} Are pairs $(X,\xxx)$ consisting of an object
$X \in \Mm$ and a morphism $\xxx: \CC\tensor{}X \to X\tensor{}\CC$
such that
\begin{eqnarray}
  (X\tensor{}\Delta) \circ \xxx  &=& (\xxx\tensor{}\CC) \circ (\CC\tensor{}\xxx)
  \circ (\Delta\tensor{}X) \label{1-cell} \\
  (X\tensor{}\varepsilon) \circ \xxx &=& \varepsilon \tensor{}X. \label{1-cell'}
\end{eqnarray}
We have the following useful lemma
\begin{lemma}\label{objetos}
For every object $X \in \Mm$, the following conditions are
equivalent
\begin{enumerate}[(i)]
\item $\CC\tensor{}X$ is a $\CC$-bicomodule with a left
$\CC$-coaction $\lambda^{\CC\tensor{}X}\,=\, \Delta\tensor{}X$;

\item there is a morphism $\xxx: \CC\tensor{}X \to
X\tensor{}\CC$ satisfying equalities \eqref{1-cell} and
\eqref{1-cell'}.
\end{enumerate}
\end{lemma}
\begin{proof}
$(ii) \Rightarrow (i)$. Take the right $\CC$-coaction
$\rho^{\CC\tensor{}X}\,=\, (\CC\tensor{}\xxx) \circ
(\Delta\tensor{}X)$.\\ $(i) \Rightarrow (ii)$. Take $\xxx\,=\,
(\varepsilon\tensor{}X\tensor{}\CC) \circ \rho^{\CC\tensor{}X}$,
where $\rho^{\CC\tensor{}X}$ is the given right coaction.
\end{proof}
In this way the morphisms in $\rR_{\CC}^c$ are defined in their
unreduced form as follows
\paragraph{\emph{Morphisms in} $\rR_{\CC}^{\,c}$:} $$ {\rm Hom}_{\rR_{\CC}^c}\lr{ (X,\xxx),\, (X',\xxx')
} \,:=\, {\rm Hom}_{\CC-\CC}\lr{ \CC\tensor{}X,\,\CC\tensor{}X'},$$
where $\CC\tensor{}X$ and $\CC\tensor{}X'$ are endowed with the
structure of $\CC$-bicomodule defined in Lemma \ref{objetos}. That
is a morphism $\alpha: (X,\xxx) \to (X',\xxx)$ in $\rR_{\CC}^c$ is
morphism $\alpha:\CC\tensor{}X \to \CC\tensor{}X'$ in $\Mm$
satisfying
\begin{eqnarray}
  (\Delta\tensor{}X') \circ \alpha &=&  (\CC\tensor{}\alpha) \circ (\Delta\tensor{}X), \label{RC2-cell}\\
  (\CC\tensor{}\xxx') \circ (\Delta\tensor{}X') \circ \alpha &=&
  (\alpha\tensor{}\CC) \circ (\CC\tensor{}\xxx) \circ
  (\Delta\tensor{}X). \label{RC2-cell'}
\end{eqnarray}

\paragraph{\emph{The multiplications of} $\rR_{\CC}^{\,c}$:} Let
$\alpha: (X,\xxx) \to (X',\xxx')$ and $\beta: (Y,\yyy) \to
(Y',\yyy')$ two morphisms in $\rR_{\CC}^c$. One can easily proves
that
\begin{eqnarray*}
  (X,\xxx) \tensor{\CC}(Y,\yyy) &:=& \lr{ X\tensor{}Y,\,
(X\tensor{}\yyy) \circ (\xxx\tensor{}Y)} \\
  (X',\xxx') \tensor{\CC}(Y',\yyy') &:=& \lr{ X'\tensor{}Y',\,
(X'\tensor{}\yyy') \circ (\xxx'\tensor{}Y')}
\end{eqnarray*}
are also an objects of the category $\rR_{\CC}^c$, which defines the
horizontal multiplication. The vertical one is defined as a composed
morphism:
\begin{eqnarray*}
  (\alpha \tensor{\CC}\beta) & := & (\CC\tensor{}X'\tensor{}\varepsilon\tensor{}Y') \circ (\CC\tensor{}X'\tensor{}\beta)
   \circ (\CC\tensor{}\xxx'\tensor{}Y) \circ (\CC\tensor{}\alpha\tensor{}Y) \circ (\Delta\tensor{}X\tensor{}Y)  \\
   &=& (\CC\tensor{}X'\tensor{}\varepsilon\tensor{}Y') \circ
   (\alpha\tensor{}\beta) \circ (\CC\tensor{}\xxx\tensor{}Y) \circ
   (\Delta\tensor{}X\tensor{}Y).
\end{eqnarray*}
Lastly, the identity object of the multiplication $-\tensor{\CC}-$
is given by the pair $(\II, \CC)$ (here $\CC$ denotes the identity
morphism of $\CC$ in $\Mm$).

\subsection{The monoidal category $\lL_{\CC}^c$}\label{LC} $ $ \\
This is the left Eilenberg-Moore category associated to $\CC$

\paragraph{\emph{Objects of}  $\lL_{\CC}^{\,c}$:} Are pairs $(\ppp,P)$ consisting of an object
$P \in \Mm$ and a morphism $\ppp: P\tensor{}\CC \to \CC\tensor{}P$
such that
\begin{eqnarray}
  (\Delta\tensor{}P) \circ \ppp  &=& (\CC\tensor{}\ppp) \circ (\ppp\tensor{}\CC)
  \circ (P\tensor{}\Delta) \label{l1-cell} \\
  (\varepsilon\tensor{}P) \circ \ppp &=& P \tensor{}\varepsilon. \label{l1-cell'}
\end{eqnarray}
As before one can easily checks the following lemma
\begin{lemma}\label{lobjetos}
For every object $P \in \Mm$, the following conditions are
equivalent
\begin{enumerate}[(i)]
\item $P\tensor{}\CC$ is a $\CC$-bicomodule with a right
$\CC$-coaction $\rho^{P\tensor{}\CC}\,=\, P\tensor{}\Delta$;

\item there is a morphism $\ppp: P\tensor{}\CC \to
\CC\tensor{}P$ satisfying equalities \eqref{l1-cell} and
\eqref{l1-cell'}.
\end{enumerate}
\end{lemma}
In this way the morphisms in $\lL_{\CC}^c$ are defined in their
unreduced form as follows
\paragraph{\emph{Morphisms in} $\lL_{\CC}^{\,c}$:} $$ {\rm Hom}_{\lL_{\CC}^c}\lr{ (\ppp,P),\, (\ppp',P')
} \,:=\, {\rm Hom}_{\CC-\CC}\lr{ P\tensor{}\CC,\,P'\tensor{}\CC},$$
where $P\tensor{}\CC$ and $P'\tensor{}\CC$ are endowed with the
structure of $\CC$-bicomodule defined in Lemma \ref{lobjetos}. That
is a morphism $\gamma: (\ppp,P) \to (\ppp',P')$ in $\lL_{\CC}^c$ is
morphism $\gamma:P\tensor{}\CC \to P'\tensor{}\CC$ in $\Mm$
satisfying
\begin{eqnarray}
  (P'\tensor{}\Delta) \circ \gamma &=&  (\gamma\tensor{}\CC) \circ (P\tensor{}\Delta), \label{LC2-cell}\\
  (\ppp'\tensor{}\CC) \circ (P'\tensor{}\Delta) \circ \gamma &=&
  (\CC\tensor{}\gamma) \circ (\ppp\tensor{}\CC) \circ
  (P\tensor{}\Delta). \label{LC2-cell'}
\end{eqnarray}

\paragraph{\emph{The multiplications of} $\lL_{\CC}^{\,c}$:} Let
$\gamma: (\ppp,P) \to (\ppp',P')$ and $\sigma: (\qqq,Q) \to
(\qqq',Q')$ two morphisms in $\lL_{\CC}^c$. One can easily proves
that
\begin{eqnarray*}
  (\ppp,P) \otimes^{\CC} (\qqq,Q) &:=& \lr{ (\ppp\tensor{}Q) \circ (P\tensor{}\qqq),\, P\tensor{}Q} \\
  (\ppp',P') \otimes^{\CC}(\qqq',Q') &:=& \lr{ (\ppp'\tensor{}Q') \circ (P'\tensor{}\qqq'),\, P'\tensor{}Q'}
\end{eqnarray*}
are also an objects of the category $\lL_{\CC}^c$, which leads to
the horizontal multiplication. The vertical one is defined as the
composed morphism:
\begin{eqnarray*}
  (\gamma \otimes^{\CC}\sigma) & := & (P'\tensor{}\varepsilon\tensor{}Q'\tensor{}\CC) \circ (\gamma\tensor{}Q'\tensor{}\CC)
   \circ (P\tensor{}\qqq'\tensor{}\CC) \circ (P\tensor{}\sigma\tensor{}\CC) \circ (P\tensor{}Q\tensor{}\Delta)  \\
   &=& (P'\tensor{}\varepsilon\tensor{}Q'\tensor{}\CC) \circ
   (\gamma\tensor{}\sigma) \circ (P\tensor{}\qqq\tensor{}\CC) \circ
   (P\tensor{}Q\tensor{}\Delta).
\end{eqnarray*}
Lastly, the identity object of the multiplication $-\tensor{\CC}-$
is given by the pair $(\CC,\II)$ (here $\CC$ denotes the identity
morphism of $\CC$ in $\Mm$).

\subsection{Cowreath and their products.}\label{Wreath}

$\CC$ still denotes a comonoid in $\Mm$, $\rR_{\CC}^c$ and
$\lL_{\CC}^c$ are the monoidal categories defined, respectively, in
subsections \ref{RC} and \ref{LC}. The notion of wreath was
introduced in \cite{Lack/Street:2002} in the general context of
$2$-categories, in the monoidal case they are defined as follows:

\begin{definition}\label{Def-Wreath}
Let $\CC$ be a comonoid in a strict monoidal category $\Mm$. A
\emph{right cowreath over} $\CC$ (or \emph{right $\CC$-cowreath}) is
a comonoid in the monoidal category $\rR_{\CC}^c$. A \emph{right
wreath over} $\CC$ (or \emph{right $\CC$-wreath}) is a monoid in the
monoidal category $\rR_{\CC}^c$. The left versions of these
definitions are obtained in the monoidal category $\lL_{\CC}^c$.
\end{definition}

The following gives, in terms of the multiplication of $\Mm$, a
simplest and equivalent definition of cowreath.

\begin{proposition}\label{Def-equ}
Let $\CC$ be a comonoid in a strict monoidal category $\Mm$ and
$(\RR,\rrr)$ an object of the category $\rR_{\CC}^c$. The following
statements are equivalent
\begin{enumerate}[(i)]
\item $(\RR,\rrr)$ is a right $\CC$-cowreath;

\item There is a $\CC$-bicomodules morphisms $\xi: \CC\tensor{}\RR \to \CC$
and $\delta: \CC\tensor{}\RR \to \CC\tensor{}\RR\tensor{}\RR$
converting the following diagrams commutative
$$
\xy *+{\CC\tensor{}\RR}="p",
p+<3cm,0pt>*+{\CC\tensor{}\RR\tensor{}\RR}="1",
p+<3cm,-2cm>*+{\CC\tensor{}\RR}="2",{"p" \ar@{=} "2"}, {"p"
\ar@{->}^-{\delta} "1"}, {"1" \ar@{->}^-{\xi\tensor{}\RR} "2"}
\endxy \qquad \xy *+{\CC\tensor{}\RR}="p",
p+<4cm,0pt>*+{\CC\tensor{}\RR\tensor{}\RR}="1",
p+<4cm,-2cm>*+{\RR\tensor{}\CC\tensor{}\RR}="2",
p+<0pt,-2cm>*+{\RR\tensor{}\CC}="3", {"p" \ar@{->}_-{\rrr} "3"},
{"p" \ar@{->}^-{\delta} "1"}, {"1" \ar@{->}^-{\rrr\tensor{}\RR}
"2"}, {"2" \ar@{->}^-{\RR\tensor{}\xi} "3"}
\endxy $$
$$
\xy *+{\CC\tensor{}\RR}="p",
p+<4cm,0pt>*+{\CC\tensor{}\RR\tensor{}\RR}="1",
p+<8cm,0pt>*+{\CC\tensor{}\RR\tensor{}\RR\tensor{}\RR}="2",
p+<0pt,-2cm>*+{\CC\tensor{}\RR\tensor{}\RR}="3",
p+<4cm,-2cm>*+{\RR\tensor{}\CC\tensor{}\RR}="4",
p+<8cm,-2cm>*+{\RR\tensor{}\CC\tensor{}\RR\tensor{}\RR}="5",  {"p"
\ar@{->}^-{\delta} "1"}, {"p" \ar@{->}_-{\delta} "3"}, {"1"
\ar@{->}^-{\delta\tensor{}\RR} "2"}, {"2"
\ar@{->}^{\rrr\tensor{}\RR\tensor{}\RR} "5"}, {"3"
\ar@{->}_-{\rrr\tensor{}\RR} "4"}, {"4"
\ar@{->}_-{\RR\tensor{}\delta} "5"}
\endxy
$$
\end{enumerate}
\end{proposition}
\begin{proof}
Analogue to that of \cite[Proposition 2.2]{Kaoutit:2007pr}.
\end{proof}

The cowreath products was introduced in \cite{Lack/Street:2002} for
comonads in a general $2$-categories. In the particular case of
strict monoidal categories this product is expressed  by the
following

\begin{proposition}\label{proct-cowreath}
Let $\CC$ be a comonoid in $\Mm$, and $(\RR,\rrr)$ a right
$\CC$-cowreath with structure morphisms $\xi: \CC\tensor{}\RR \to
\CC$ and $\delta: \CC\tensor{}\RR \to \CC\tensor{}\RR\tensor{}\RR$.
The object $\CC\tensor{}\RR$ admits a structure of comonoid with
comultiplication and counit given by
$$\xymatrix@R=10pt@C=40pt{\bd{\Delta}: \CC\tensor{}\RR \ar@{->}^-{\Delta\tensor{}\RR}[r] &
\CC\tensor{}\CC\tensor{}\RR \ar@{->}^-{\CC\tensor{}\delta}[r] &
\CC\tensor{}\CC\tensor{}\RR\tensor{}\RR
\ar@{->}^-{\CC\tensor{}\rrr\tensor{}\RR}[r] &
\CC\tensor{}\RR\tensor{}\CC\tensor{}\RR, \\ \bd{\varepsilon}:
\CC\tensor{}\RR \ar@{->}^-{\xi}[r] & \CC \ar@{->}^-{\varepsilon}[r]
& \II. & }$$ Moreover, with this comonoid structure the morphism
$\xi: \CC\tensor{}\RR \to \CC$ becomes a morphism of comonoid.
\end{proposition}
\begin{proof}
Straightforward.
\end{proof}
The comonoid $\CC\tensor{}\RR$ of the previous Proposition is
referred as \emph{the cowreath product} of $\CC$ by $\RR$.

\begin{remark}\label{DL-c}
Notice that the object $\RR$ occurring in Proposition \ref{Def-equ}
need not be a comonoid. However, if $\RR$ is it self a comonoid with
a structure morphisms $\Delta': \RR \to \RR\tensor{}\RR$,
$\varepsilon': \RR \to \II$  such that the pair  $(\rrr,\CC)$
belongs to the left monoidal category $\lL_{\RR}^c$, then the
morphisms $\xi:=\CC\tensor{}\varepsilon'$ and
$\delta:=\CC\tensor{}\Delta'$ endow $(\RR,\rrr)$ with a structure of
right $\CC$-cowreath while $\xi':=\varepsilon\tensor{}\RR$, and
$\delta':=\Delta\tensor{}\RR$ gave to $(\rrr,\CC)$ a structure of
left $\RR$-cowreath. Furthermore, by Proposition
\ref{proct-cowreath}, the morphisms $\xi: \CC\tensor{}\RR \to \CC$
and $\xi':\CC\tensor{}\RR \to \RR$ are in fact a comonoids
morphisms.\\ In this way the morphism $\rrr: \CC\tensor{}\RR \to
\RR\tensor{}\CC$ should satisfies the following equalities
\begin{eqnarray}
  (\RR\tensor{}\Delta) \circ \rrr &=& (\rrr\tensor{}\CC) \circ (\CC\tensor{}\rrr)
  \circ (\Delta\tensor{}\RR) \label{CDL-1}  \\
   (\RR\tensor{}\varepsilon) \circ \rrr &=& \varepsilon\tensor{}\RR \label{CDL-2} \\
  (\Delta'\tensor{}\CC) \circ \rrr &=& (\RR\tensor{}\rrr) \circ
  (\rrr\tensor{}\RR) \circ (\CC\tensor{}\Delta') \label{CDL-3}  \\
   (\varepsilon'\tensor{}\CC) \circ \rrr &=&
   \CC\tensor{}\varepsilon'. \label{CDL-4}
\end{eqnarray}
A morphism satisfying the four previous equalities is called a
\emph{comonoid distributive law} from $\CC$ to $\RR$, see
\cite{Beck:1969} for the original definition.
\end{remark}

A cowreath product satisfies a universal property in the following
sense

\begin{proposition}\label{Univ-coPrp}
Let $(\CC,\Delta,\varepsilon)$ be a comonoid in a strict monoidal
category $\Mm$, and $(\RR,\rrr)$ a $\CC$-cowreath with a structure
morphisms $\xi: \CC\tensor{}\RR \to \CC$ and $\delta:\CC\tensor{}\RR
\to \CC\tensor{}\RR\tensor{}\RR$.

Let $(\DD,\Delta',\varepsilon')$ be a comonoid with a comonoid
morphism $\alpha:\DD \to \CC$ and with morphism $\beta: \DD \to \RR$
satisfying
\begin{eqnarray}
  \xi \circ (\CC\tensor{}\beta) &=& \CC\tensor{}\varepsilon' \label{Hy-1}\\
  \delta \circ (\CC\tensor{}\beta) &=&
  (\CC\tensor{}\beta\tensor{}\beta) \circ (\CC\tensor{}\Delta') \label{Hy-2}
\end{eqnarray}
Assume that $\alpha$ and $\beta$ satisfy the equality
\begin{eqnarray}
  \rrr \circ (\alpha\tensor{}\beta) \circ \Delta' &=& (\beta\tensor{}\alpha)
  \circ \Delta' \label{Hy-3}
\end{eqnarray}
then there exists a unique comonoid morphism $\gamma: \DD \to
\CC\tensor{}\RR$ such that $\xi \circ \gamma\,=\, \alpha$ and
$(\varepsilon\tensor{}\RR) \circ \gamma\,=\,\beta$.
\end{proposition}
\begin{proof}
If there exists such a morphism, then it should be unique by the
following computations
\begin{eqnarray*}
  (\alpha\tensor{}\beta) \circ \Delta'
   &=& \lr{(\xi \circ \gamma)\tensor{}\lr{(\varepsilon\tensor{}\RR) \circ \gamma}} \circ \Delta' \\
   &=& (\xi\tensor{}\varepsilon\tensor{}\RR) \circ (\gamma\tensor{}\gamma) \circ \Delta' \\
   &=& (\xi\tensor{}\varepsilon\tensor{}\RR) \circ (\CC\tensor{}\rrr\tensor{}\RR) \circ
   (\CC\tensor{}\delta) \circ (\Delta\tensor{}\RR) \circ \gamma \\
   &=& (\xi\tensor{}\RR) \circ (\CC\tensor{}\RR\tensor{}\varepsilon\tensor{}\RR) \circ (\CC\tensor{}\rrr\tensor{}\RR) \circ
   (\CC\tensor{}\delta) \circ (\Delta\tensor{}\RR) \circ \gamma \\
   &\overset{\eqref{1-cell'}}{=}& (\xi\tensor{}\RR) \circ (\CC\tensor{}\varepsilon\tensor{}\RR\tensor{}\RR) \circ
   (\CC\tensor{}\delta) \circ (\Delta\tensor{}\RR) \circ \gamma \\
   &\overset{\eqref{RC2-cell}}{=}& (\xi\tensor{}\RR) \circ (\CC\tensor{}\varepsilon\tensor{}\RR\tensor{}\RR) \circ
   (\Delta\tensor{}\RR\tensor{}\RR) \circ \delta \circ \gamma \\
   &=&  (\xi\tensor{}\RR) \circ \delta \circ \gamma \,\,
   \overset{\ref{Def-equ}}{=}\,\, \gamma.
\end{eqnarray*}
Since by hypothesis $\xi \circ (\alpha\tensor{}\beta) \circ
\Delta'\,=\, \alpha$ and $(\varepsilon\tensor{}\RR) \circ
(\alpha\tensor{}\beta) \circ \Delta'\,=\, \beta$, it suffice to show
that $(\alpha\tensor{}\beta) \circ \Delta'$ is a comonoid morphism.
The counitary property comes out as
\begin{eqnarray*}
  \varepsilon \circ \xi \circ (\alpha\tensor{}\beta) \circ \Delta' &=&
  \varepsilon \circ \xi \circ (\CC\tensor{}\beta) \circ (\alpha\tensor{}\DD) \circ \Delta' \\
   &\overset{\eqref{Hy-1}}{=}& \varepsilon \circ (\CC\tensor{}\varepsilon') \circ (\alpha\tensor{}\DD) \circ \Delta' \\
   &=& \varepsilon \circ \alpha \circ (\DD\tensor{}\varepsilon') \circ \Delta'
   \,\,=\,\, \varepsilon'.
\end{eqnarray*}
Now the coassociativity property is obtained by the following
computations
\begin{eqnarray*}
   \bd{\Delta} \circ (\alpha\tensor{}\beta) \circ \Delta' &=&
   (\CC\tensor{}\rrr\tensor{}\RR) \circ (\CC\tensor{}\delta) \circ
   (\Delta\tensor{}\RR) \circ (\alpha\tensor{}\beta) \circ \Delta'  \\
   &\overset{\eqref{RC2-cell}}{=}& (\CC\tensor{}\rrr\tensor{}\RR) \circ (\Delta\tensor{}\RR\tensor{}\RR) \circ
   \delta \circ (\alpha\tensor{}\beta) \circ \Delta' \\
   &=& (\CC\tensor{}\rrr\tensor{}\RR) \circ (\Delta\tensor{}\RR\tensor{}\RR) \circ
   \delta \circ (\CC\tensor{}\beta) \circ (\alpha\tensor{}\DD) \circ \Delta'  \\
   &\overset{\eqref{Hy-2}}{=}& (\CC\tensor{}\rrr\tensor{}\RR) \circ (\Delta\tensor{}\RR\tensor{}\RR) \circ
   (\CC\tensor{}\beta\tensor{}\beta) \circ (\CC\tensor{}\Delta') \circ (\alpha\tensor{}\DD) \circ \Delta'  \\
   &=& (\CC\tensor{}\rrr\tensor{}\RR) \circ (\Delta\tensor{}\RR\tensor{}\RR) \circ
   (\alpha\tensor{}\RR\tensor{}\RR)\circ (\DD\tensor{}\beta\tensor{}\beta)  \circ (\DD\tensor{}\Delta') \circ \Delta' \\
   &=& (\CC\tensor{}\rrr\tensor{}\RR) \circ (\alpha\tensor{}\alpha\tensor{}\RR\tensor{}\RR) \circ
   (\Delta'\tensor{}\RR\tensor{}\RR)\circ (\DD\tensor{}\beta\tensor{}\beta)  \circ (\DD\tensor{}\Delta') \circ \Delta' \\
   &=& (\alpha\tensor{}\RR\tensor{}\CC\tensor{}\RR)
   \circ \lr{\DD\tensor{}(\rrr \circ (\alpha\tensor{}\beta))\tensor{}\RR} \circ (\Delta'\tensor{}\DD\tensor{}\RR)
   \circ (\DD\tensor{}\DD\tensor{}\beta) \circ (\Delta'\tensor{}\DD) \circ \Delta' \\
   &=& (\alpha\tensor{}\RR\tensor{}\CC\tensor{}\RR)
   \circ \lr{\DD\tensor{}(\rrr \circ (\alpha\tensor{}\beta))\tensor{}\RR} \circ (\Delta'\tensor{}\DD\tensor{}\RR)
   \circ (\Delta'\tensor{}\RR) \circ (\DD\tensor{}\beta) \circ \Delta' \\
   &=& (\alpha\tensor{}\RR\tensor{}\CC\tensor{}\RR)
   \circ \lr{\DD\tensor{}\lr{\rrr \circ (\alpha\tensor{}\beta) \circ \Delta'}\tensor{}\RR}
   \circ (\Delta'\tensor{}\RR) \circ (\DD\tensor{}\beta) \circ \Delta' \\
   &\overset{\eqref{Hy-3}}{=}& (\alpha\tensor{}\RR\tensor{}\CC\tensor{}\RR)
   \circ (\DD\tensor{}\beta\tensor{}\alpha\tensor{}\RR) \circ
   (\DD\tensor{}\Delta'\tensor{}\RR) \circ (\Delta'\tensor{}\RR) \circ (\DD\tensor{}\beta) \circ \Delta'  \\
   &=& (\alpha\tensor{}\beta\tensor{}\alpha\tensor{}\RR)
   \circ (\DD\tensor{}\Delta'\tensor{}\RR) \circ
   (\DD\tensor{}\DD\tensor{}\beta) \circ (\Delta'\tensor{}\DD) \circ \Delta' \\
   &=& (\alpha\tensor{}\beta\tensor{}\alpha\tensor{}\beta)
   \circ (\DD\tensor{}\Delta'\tensor{}\DD) \circ (\Delta'\tensor{}\DD) \circ \Delta' \\
   &=& (\alpha\tensor{}\beta\tensor{}\alpha\tensor{}\beta)
   \circ (\Delta'\tensor{}\Delta')  \circ \Delta'\,\, =\,\, \lr{
   \lr{ (\alpha\tensor{}\beta) \circ \Delta'} \tensor{} \lr{(\alpha\tensor{}\beta) \circ
   \Delta'}} \circ \Delta'.
\end{eqnarray*}
\end{proof}

In what follows we announce the analogue notion for a given monoid
in a strict monoidal category. So consider a monoid $(\AA,\mu,\eta)$
in $\Mm$. We start by defining the right and left Eilenberg-Moore
monoidal categories attached to $\AA$

\subsection{The monoidal category $\rR_{\AA}^{\,a}$}\label{RA} $ $ \\
This is the right Eilenberg-Moore monoidal category associated to
the monoid $\AA$ (see \cite{Lack/Street:2002}), and defined as
follows

\paragraph{\emph{Objects of}  $\rR_{\AA}^{\,a}$:} Are pairs $(U,\uuu)$ consisting of an object
$U \in \Mm$ and a morphism $\uuu: \AA\tensor{}U \to U\tensor{}\AA$
such that
\begin{eqnarray}
  \uuu \circ (\mu\tensor{}U)  &=& (U\tensor{}\mu) \circ (\uuu\tensor{}\AA)
  \circ (\AA\tensor{}\uuu) \label{RA1-cell} \\
  \uuu \circ (\eta\tensor{}U)  &=& U\tensor{}\eta. \label{RA1-cell'}
\end{eqnarray}
We have the following lemma
\begin{lemma}\label{RA-objetos}
For every object $U \in \Mm$, the following conditions are
equivalent
\begin{enumerate}[(i)]
\item $U\tensor{}\AA$ is an $\AA$-bimodule with a right
$\AA$-action $\Sf{r}_{U\tensor{}\AA}\,=\, U\tensor{}\mu$;

\item there is a morphism $\uuu: \AA\tensor{}U \to
U\tensor{}\AA$ satisfying equalities \eqref{RA1-cell} and
\eqref{RA1-cell'}.
\end{enumerate}
\end{lemma}
\begin{proof}
$(ii) \Rightarrow (i)$. Take the left $\AA$-action
$\Sf{l}_{U\tensor{}\AA}\,=\, (U\tensor{}\mu) \circ (\uuu\tensor{}\AA)$.\\
$(i) \Rightarrow (ii)$. Take $\uuu\,=\, \Sf{l}_{U\tensor{}\AA} \circ
(\AA\tensor{}U\tensor{}\eta)$, where $\Sf{l}_{U\tensor{}\AA}$ is the
given left action.
\end{proof}
In this way the morphisms in $\rR_{\AA}^a$ are defined in their
unreduced form as follows
\paragraph{\emph{Morphisms in} $\rR_{\AA}^{\,a}$:} $$ {\rm Hom}_{\rR_{\AA}^a}\lr{ (U,\uuu),\, (U',\uuu')
} \,:=\, {\rm Hom}_{\AA-\AA}\lr{ U\tensor{}\AA,\,U'\tensor{}\AA},$$
where $U\tensor{}\AA$ and $U'\tensor{}\AA$ are endowed with the
structure of $\AA$-bimodule defined in Lemma \ref{RA-objetos}. That
is amorphism $\nu: (U,\uuu) \to (U',\uuu')$ in $\rR^a_{\AA}$ is a
morphism $\nu : U\tensor{}\AA \to U'\tensor{}\AA$ satisfying
\begin{eqnarray}
  \nu \circ (U\tensor{}\mu) &=&  (U'\tensor{}\mu) \circ (\nu\tensor{}\AA), \label{RA2-cell}\\
  \nu \circ (U\tensor{}\mu) \circ (\uuu\tensor{}\AA) &=&
  (U'\tensor{}\mu) \circ (\uuu'\tensor{}\AA) \circ
  (\AA\tensor{}\nu). \label{RA2-cell'}
\end{eqnarray}

\paragraph{\emph{The multiplications of} $\rR_{\AA}^{\,a}$:} Let
$\nu: (U,\uuu) \to (U',\uuu')$ and $\upsilon: (V,\vvv) \to
(V',\vvv')$ two morphisms in $\rR_{\AA}^a$. One can easily proves
that
\begin{eqnarray*}
  (U,\uuu) \tensor{\AA}(V,\vvv) &:=& \lr{ U\tensor{}V,\,
(U\tensor{}\vvv) \circ (\uuu\tensor{}V)} \\
  (U',\uuu') \tensor{\AA}(V',\vvv') &:=& \lr{ U'\tensor{}V',\,
(U'\tensor{}\vvv') \circ (\uuu'\tensor{}V')}
\end{eqnarray*}
are also an objects of the category $\rR_{\AA}^a$, which gives the
horizontal multiplication. The vertical one is defined by the
composition:
\begin{eqnarray*}
  (\nu \tensor{\AA}\upsilon) & := & (U'\tensor{}V'\tensor{}\mu) \circ (U'\tensor{}\upsilon\tensor{}\AA)
   \circ (U'\tensor{}\vvv\tensor{}\AA) \circ (\nu\tensor{}V\tensor{}\AA) \circ (U\tensor{}\eta\tensor{}V\tensor{}\AA)  \\
   &=& (U'\tensor{}V'\tensor{}\mu) \circ
   (U'\tensor{}\vvv'\tensor{}\AA) \circ (\nu\tensor{}\upsilon) \circ
   (U\tensor{}\eta\tensor{}V\tensor{}\AA)
\end{eqnarray*}
Lastly, the identity object of the multiplication $-\tensor{\AA}-$
is given by the pair $(\II, \AA)$ (here $\AA$ denotes the identity
morphism of $\AA$ in $\Mm$).

\subsection{The monoidal category $\lL_{\AA}^a$}\label{LA} $ $ \\
This is the left Eilenberg-Moore category associated to $\AA$, and
defined as follows

\paragraph{\emph{Objects of}  $\lL_{\AA}^{\,a}$:} Are pairs $(\mmm,M)$ consisting of an object
$M \in \Mm$ and a morphism $\mmm: M\tensor{}\AA \to \AA\tensor{}M$
such that
\begin{eqnarray}
  \mmm \circ (M\tensor{}\mu)  &=& (\mu\tensor{}M) \circ
  (\AA\tensor{}\mmm) \circ (\mmm\tensor{}\AA) \label{LA1-cell} \\
  \mmm \circ (M\tensor{}\eta)  &=& \eta \tensor{}M. \label{LA1-cell'}
\end{eqnarray}
As before one can easily check the following lemma
\begin{lemma}\label{LA-objetos}
For every object $M \in \Mm$, the following conditions are
equivalent
\begin{enumerate}[(i)]
\item $\AA\tensor{}M$ is an $\AA$-bimodule with a left
$\AA$-action $\Sf{l}_{\AA\tensor{}M}\,=\, \mu\tensor{}M$;

\item there is a morphism $\mmm: M\tensor{}\AA \to
\AA\tensor{}M$ satisfying equalities \eqref{LA1-cell} and
\eqref{LA1-cell'}.
\end{enumerate}
\end{lemma}
In this way the morphisms in $\lL_{\AA}^a$ are defined in their
unreduced form as follows
\paragraph{\emph{Morphisms in} $\lL_{\AA}^{\,a}$:} $$ {\rm Hom}_{\lL_{\AA}^a}\lr{ (\mmm,M),\, (\mmm',M')
} \,:=\, {\rm Hom}_{\AA-\AA}\lr{ \AA\tensor{}M,\,\AA\tensor{}M'},$$
where $\AA\tensor{}M$ and $\AA\tensor{}M'$ are endowed with the
structure of $\AA$-bimodule defined in Lemma \ref{LA-objetos}. That
is amorphism $\theta: (\mmm,M) \to (\mmm',M')$ in $\lL^a_{\AA}$ is a
morphism $\theta : \AA\tensor{}M \to \AA\tensor{}M'$ satisfying
\begin{eqnarray}
  \theta \circ (\mu\tensor{}M) &=&  (\mu\tensor{}M') \circ (\AA\tensor{}\theta), \label{LA2-cell}\\
  \theta \circ (\mu\tensor{}M) \circ (\AA\tensor{}\mmm) &=&
  (\mu\tensor{}M') \circ (\AA\tensor{}\mmm') \circ
  (\theta\tensor{}\AA). \label{LA2-cell'}
\end{eqnarray}

\paragraph{\emph{The multiplications of} $\lL_{\AA}^{\,a}$:} Let
$\theta: (\mmm,M) \to (\mmm',M')$ and $\vartheta: (\nnn,N) \to
(\nnn',N')$ two morphisms in $\lL_{\AA}^a$. One can easily prove
that
\begin{eqnarray*}
  (\mmm,M) \otimes^{\AA} (\nnn,N) &:=& \lr{ (\mmm\tensor{}N) \circ (M\tensor{}\nnn),\, M\tensor{}N} \\
  (\mmm',M') \otimes^{\AA}(\nnn',N') &:=& \lr{ (\nnn'\tensor{}N') \circ (M'\tensor{}\nnn'),\, M'\tensor{}N'}
\end{eqnarray*}
are also an objects of the category $\lL_{\AA}^a$, which leads to
the horizontal multiplication. The vertical one is defined as the
composed morphism:
\begin{eqnarray*}
  (\theta \otimes^{\AA}\vartheta) & := & (\mu\tensor{}M'\tensor{}N') \circ (\AA\tensor{}\theta\tensor{}N')
   \circ (\AA\tensor{}\mmm\tensor{}N') \circ (\AA\tensor{}M\tensor{}\vartheta)
   \circ (\AA\tensor{}M\tensor{}\eta\tensor{}N)  \\
   &=& (\mu\tensor{}M'\tensor{}N')  \circ
   (\AA\tensor{}\mmm'\tensor{}N') \circ (\theta\tensor{}\vartheta) \circ
   (\AA\tensor{}M\tensor{}\eta\tensor{}N).
\end{eqnarray*}
Lastly, the identity object of the multiplication $-\otimes^{\AA}-$
is given by the pair $(\AA,\II)$.

\subsection{Wreaths and their products.}\label{A-Wreath}

$\AA$ still denotes a monoid in $\Mm$, $\rR_{\AA}^a$ and
$\lL_{\AA}^a$ are the monoidal categories defined, respectively, in
subsections \ref{RA} and \ref{LA}. As in the comonoid case, we have

\begin{definition}\label{A-Def-Wreath}
Let $\AA$ be a monoid in a strict monoidal category $\Mm$. A
\emph{right wreath over} $\AA$ (or \emph{right $\AA$-wreath}) is a
monoid in the monoidal category $\rR_{\AA}^a$. A \emph{right
cowreath over} $\AA$ (or \emph{right $\AA$-cowreath}) is a comonoid
in the monoidal category $\rR_{\AA}^a$. The left versions of these
notions are defined in the monoidal category $\lL_{\CC}^c$.
\end{definition}

The following gives, in terms of the multiplication of $\Mm$, a
simplest and equivalent definition of wreath.

\begin{proposition}\label{A-Def-equ}
Let $\AA$ be a monoid in a strict monoidal category $\Mm$ and
$(\TT,\ttt)$ an object of the category $\rR_{\AA}^a$. The following
statements are equivalent
\begin{enumerate}[(i)]
\item $(\TT,\ttt)$ is a right $\AA$-wreath;

\item There is an $\AA$-bimodules morphisms $\zeta: \AA \to \TT\tensor{}\AA$
and $\nu: \TT\tensor{}\TT\tensor{}\AA \to \TT\tensor{}\AA$ rendering
the following diagrams commutative
$$
\xy *+{\TT\tensor{}\TT\tensor{}\AA}="p",
p+<3cm,0pt>*+{\TT\tensor{}\AA}="1",
p+<0pt,-2cm>*+{\TT\tensor{}\AA}="2",{"2"
\ar@{->}^-{\TT\tensor{}\zeta} "p"}, {"p" \ar@{->}^-{\nu} "1"}, {"1"
\ar@{=} "2"}
\endxy \qquad \xy *+{\TT\tensor{}\TT\tensor{}\AA}="p",
p+<4cm,0pt>*+{\TT\tensor{}\AA}="1",
p+<4cm,-2cm>*+{\AA\tensor{}\TT}="2",
p+<0pt,-2cm>*+{\TT\tensor{}\AA\tensor{}\TT}="3", {"3"
\ar@{->}^-{\TT\tensor{}\ttt} "p"}, {"p" \ar@{->}^-{\nu} "1"}, {"2"
\ar@{->}_-{\ttt} "1"}, {"2" \ar@{->}^-{\zeta\tensor{}\TT} "3"}
\endxy $$
$$
\xy *+{\TT\tensor{}\TT\tensor{}\AA\tensor{}\TT}="p",
p+<4cm,0pt>*+{\TT\tensor{}\AA\tensor{}\TT}="1",
p+<8cm,0pt>*+{\TT\tensor{}\TT\tensor{}\AA}="2",
p+<0pt,-2cm>*+{\TT\tensor{}\TT\tensor{}\TT\tensor{}\AA}="3",
p+<4cm,-2cm>*+{\TT\tensor{}\TT\tensor{}\AA}="4",
p+<8cm,-2cm>*+{\TT\tensor{}\AA}="5",  {"p"
\ar@{->}^-{\nu\tensor{}\TT} "1"}, {"p"
\ar@{->}_-{\TT\tensor{}\TT\tensor{}\ttt} "3"}, {"1"
\ar@{->}^-{\TT\tensor{}\ttt} "2"}, {"2" \ar@{->}^{\nu} "5"}, {"3"
\ar@{->}_-{\TT\tensor{}\nu} "4"}, {"4" \ar@{->}_-{\nu} "5"}
\endxy
$$
\end{enumerate}
\end{proposition}
\begin{proof}
Is left to the reader.
\end{proof}

The wreath products is expressed  by the following

\begin{proposition}\label{product-wreath}
Let $(\AA,\mu,\eta)$ be a monoid in $\Mm$, and $(\TT,\ttt)$ a right
$\AA$-wreath with structure morphisms $\zeta: \AA \to
\TT\tensor{}\AA$ and $\nu: \TT\tensor{}\TT\tensor{}\AA \to
\TT\tensor{}\AA$. The object $\TT\tensor{}\AA$ admits a structure of
monoid with multiplication and unit given by
$$\xymatrix@R=10pt@C=40pt{\mu': \TT\tensor{}\AA\tensor{}\TT\tensor{}\AA \ar@{->}^-{\TT\tensor{}\ttt\tensor{}\AA}[r] &
\TT\tensor{}\TT\tensor{}\AA\tensor{}\AA
\ar@{->}^-{\nu\tensor{}\AA}[r] & \TT\tensor{}\AA\tensor{}\AA
\ar@{->}^-{\TT\tensor{}\mu}[r] & \TT\tensor{}\AA, \\
\eta': \II \ar@{->}^-{\eta}[r] & \AA \ar@{->}^-{\zeta}[r] &
\TT\tensor{}\AA. & }$$ Moreover, with this monoid structure the
morphism $\zeta: \AA \to \TT\tensor{}\AA$ becomes a morphism of
monoides.
\end{proposition}
\begin{proof}
Straightforward.
\end{proof}
The monoid $\TT\tensor{}\AA$ of the previous Proposition is referred
as \emph{the wreath product} of $\AA$ by $\TT$.

\begin{remark}\label{DL-A}
Notice here also that the object $\TT$ occurring in Proposition
\ref{A-Def-equ} need not be a monoid. However, if $\TT$ is it self a
monoid with a structure morphisms $\mu': \TT\tensor{}\TT \to \TT$,
$\eta': \II \to \TT$  such that the pair  $(\ttt,\AA)$ belongs to
the left monoidal category $\lL_{\TT}^a$, then the morphisms
$\zeta:=\eta'\tensor{}\AA$ and $\nu:=\mu'\tensor{}\AA$ endow
$(\TT,\ttt)$ with a structure of right $\AA$-wreath while
$\zeta':=\eta\tensor{}\TT$, and $\nu':=\mu\tensor{}\AA$ gave to
$(\ttt,\AA)$ a structure of left $\TT$-wreath. Furthermore, by
Proposition \ref{product-wreath}, the morphisms $\zeta: \AA \to
\TT\tensor{}\AA$ and $\zeta':\TT \to \TT\tensor{}\AA$ are in fact a
monoids morphisms.\\ In this way the morphism $\ttt: \TT\tensor{}\AA
\to \AA\tensor{}\TT$ should satisfies the following equalities
\begin{eqnarray}
  \ttt \circ (\mu\tensor{}\TT) &=& (\TT\tensor{}\mu) \circ (\ttt\tensor{}\AA)
  \circ (\AA\tensor{}\ttt) \label{ADL-1}  \\
   \ttt \circ (\eta\tensor{}\TT) &=& \TT\tensor{}\eta \label{ADL-2} \\
  \ttt \circ (\AA\tensor{}\mu') &=& (\mu'\tensor{}\AA) \circ
  (\TT\tensor{}\ttt) \circ (\ttt\tensor{}\TT) \label{ADL-3}  \\
   \ttt \circ (\AA\tensor{}\eta') &=&
   \eta'\tensor{}\AA. \label{ADL-4}
\end{eqnarray}
A morphism satisfying the four previous equalities is called a
\emph{monoid distributive law} from $\AA$ to $\TT$, see also
\cite{Beck:1969} for the original definition.
\end{remark}

The universal property of wreath products is expressed as follows

\begin{proposition}\label{Univ-Prp}
Let $(\AA,\mu,\eta)$ be a monoid in a strict monoidal category
$\Mm$, and $(\TT,\ttt)$ a right $\AA$-wreath with a structure
morphisms $\zeta: \AA \to \TT\tensor{}\AA$,
$\nu:\TT\tensor{}\TT\tensor{}\AA \to \TT\tensor{}\AA$.

Let $(\LL,\mu',\eta')$ be a monoid with a monoid morphism
$\varphi:\AA \to \LL$ and with morphism $\psi: \TT \to \LL$
satisfying
\begin{eqnarray}
  (\psi\tensor{}\AA) \circ \zeta  &=& \eta'\tensor{}\AA \label{Hy-1A}\\
  (\psi\tensor{}\AA) \circ \nu  &=& (\mu'\tensor{}\AA) \circ (\psi\tensor{}\psi\tensor{}\AA) \label{Hy-2A}
\end{eqnarray}
Assume that $\varphi$ and $\psi$ satisfy the equality
\begin{eqnarray}
  \mu' \circ (\varphi\tensor{}\psi) &=& \mu' \circ (\psi\tensor{}\alpha)
  \circ \ttt \label{Hy-3A}
\end{eqnarray}
then there exists a unique monoid morphism $\phi: \TT\tensor{}\AA
\to \LL$ such that $\phi \circ \zeta\,=\, \varphi$ and $\phi \circ
(\TT\tensor{}\eta)\,=\,\psi$.
\end{proposition}

\section{Bi-monoid in general monoidal category.}\label{BiMon}

The letter $\Mm$ sill denotes a strict monoidal category with
multiplication $-\tensor{}-$ and identity object $\II$. Consider
$\BB$ an object of $\Mm$ such that $(\BB,\Delta,\varepsilon)$ is a
comonoid in $\Mm$ and $(\BB,\mu,\eta)$ is also a monoid in $\Mm$.
Assume that there is a morphism $\hbar: \BB\tensor{}\BB \to
\BB\tensor{}\BB$ which satisfies the following equalities:
\begin{eqnarray}
 \hbar \circ (\eta\tensor{}\BB) &=& \BB\tensor{}\eta  \label{A-1} \\
 \hbar \circ (\mu\tensor{}\BB) &=& (\BB\tensor{}\mu) \circ
  (\hbar\tensor{}\BB) \circ (\BB\tensor{}\hbar)  \label{A-2} \\
 \hbar \circ (\BB\tensor{}\eta) &=& \eta\tensor{}\BB  \label{A-3} \\
  \hbar \circ (\BB\tensor{}\mu) &=& (\mu\tensor{}\BB) \circ
  (\BB\tensor{}\hbar) \circ (\hbar\tensor{}\BB)  \label{A-4} \\
 (\BB\tensor{}\varepsilon) \circ \hbar&=& \varepsilon\tensor{}\BB  \label{C-1} \\
(\BB\tensor{}\Delta) \circ \hbar &=& (\hbar\tensor{}\BB) \circ (\BB\tensor{}\hbar) \circ (\Delta\tensor{}\BB)   \label{C-2} \\
(\varepsilon\tensor{}\BB) \circ \hbar &=& \BB\tensor{}\varepsilon   \label{C-3} \\
 (\Delta\tensor{}\BB) \circ \hbar &=& (\BB\tensor{}\hbar) \circ
  (\hbar\tensor{}\BB) \circ (\BB\tensor{}\Delta)  \label{C-4}.
\end{eqnarray}
Observe that the equalities \eqref{A-1}-\eqref{A-2} means that
$(\BB,\hbar) \in \rR_{\BB}^a$ and \eqref{A-2}-\eqref{A-3} means that
$(\hbar,\BB) \in \lL_{\BB}^a$ while $(\BB, \hbar) \in \rR_{\BB}^c$
by equalities \eqref{C-1}-\eqref{C-2} and $(\hbar,\BB) \in
\lL_{\BB}^c$ by equalities \eqref{C-2}-\eqref{C-4}. Moreover
\eqref{A-1}-\eqref{A-4} say that $\hbar$ is a monoid distributive
law from $\BB$ to $\BB$, and \eqref{C-1}-\eqref{C-4} say that
$\hbar$ is a comonoid distributive law form $\BB$ to $\BB$.

\begin{definition}\label{doble-Dl}
The morphism $\hbar: \BB\tensor{}\BB \to \BB\tensor{}\BB$ satisfying
equalities \eqref{A-1}-\eqref{C-4}, is called a \emph{double
distributive law} between the monoid $(\BB,\mu,\eta)$ and the
comonoid $(\BB,\Delta,\varepsilon)$.
\end{definition}

Notice that if $\hbar:\BB\tensor{}\BB \to \BB\tensor{}\BB$ is a
double distributive law, then $\lr{\BB\tensor{}\BB,
(\BB\tensor{}\hbar\tensor{}\BB) \circ (\Delta\tensor{}\Delta),
\varepsilon \circ (\BB\tensor{}\varepsilon)}$ is by Remark
\ref{DL-c} and Proposition \ref{proct-cowreath} a comonoid in $\Mm$,
and $\lr{\BB\tensor{}\BB, (\mu\tensor{}\mu) \circ
(\BB\tensor{}\hbar\tensor{}\BB), (\eta\tensor{}\BB) \circ \eta}$ is
by Remark \ref{DL-A} and Proposition \ref{product-wreath} a monoid
in $\Mm$. Using the both structures we can now enounce our main
result

\begin{proposition}\label{Bimonoid}
Let $\Mm$ be a strict monoidal category with multiplication
$-\tensor{}-$ and identity object $\II$. Consider a $6$-tuple
$(\BB,\Delta,\varepsilon,\mu,\eta,\hbar)$ where
$(\BB,\Delta,\varepsilon)$ is a comonoid in $\Mm$, $(\BB,\mu,\eta)$
is a monoid in $\Mm$ and $\hbar:\BB\tensor{}\BB \to \BB\tensor{}\BB$
is a double distributive law between them (i.e., satisfies
equalities \eqref{A-1}-\eqref{C-4}). The following statements are
equivalent
\begin{enumerate}[(i)]
\item $\Delta$ and $\varepsilon$ are morphisms of monoids;

\item $\mu$ and $\eta$ are morphisms of comonoids;

\item $\Delta$, $\varepsilon$, $\mu$ and $\eta$ satisfy:
\begin{flushleft}
${\rm (a)} \qquad \Delta \circ \eta \,\,=\,\, \eta \tensor{}\eta$  \\
${\rm (b)} \qquad (\mu\tensor{}\mu) \circ (\BB\tensor{}\hbar\tensor{}\BB)
\circ (\Delta\tensor{}\Delta) \,\,=\,\, \Delta \circ \mu$   \\
${\rm (c)} \qquad \varepsilon \circ \eta \,\,=\,\, \II$     \\
${\rm (d)} \qquad \varepsilon \circ \mu \,\,=\,\,
\varepsilon\tensor{}\varepsilon.$
\end{flushleft}
\end{enumerate}
\end{proposition}
\begin{proof}
It is clear from definitions that $\Delta$ is monoid morphism if and
only if the equalities $(iii)(a)-(b)$ are satisfied, and that
$\varepsilon$ is a monoid morphism if and only if the equalities
$(iii)(c)-(d)$ are verified. This shows that $(iii) \Leftrightarrow
(i)$. On the other hand $(iii)(b)-(d)$ is equivalent to say that
$\mu$ is a comonoid morphisms, and $(iii)(a)-(c)$ is equivalent to
say that $\eta$ is a comonoid morphism. This leads to the
equivalence $(iii) \Leftrightarrow (ii)$.
\end{proof}

\begin{definition}\label{bi-mono}
Let $\Mm$ be a strict monoidal category with multiplication
$-\tensor{}-$ and identity object $\II$. A \emph{bi-monoid} is a
$6$-tuple $(\BB,\Delta,\varepsilon,\mu,\eta,\hbar)$ satisfying the
equivalent conditions of Proposition \ref{Bimonoid}.
\end{definition}

\begin{remark}\label{JS}
The results stated in Sections \ref{ELM} and \ref{BiMon} can be
extended to the case of not necessary strict monoidal category by
using the multiplicative equivalence between any monoidal category
and a strict one, see \cite[Corollary 1.4]{Joyal/Street:1993}.
\end{remark}

\section{Some applications.}\label{aplication}
In what follows $\Bbbk$ denotes a commutative ring with $1$, and
$\rmod{\Bbbk}$ its category of modules.

\subsection{Compatibility condition between ring and coring.}
Let $R$ be $\Bbbk$-algebra, all bimodules are assumed to central
$\Bbbk$-bimodules, and their category will be denoted by
$\bimod{R}$. This is a monoidal category with multiplication the two
variables functor $-\tensor{R}-$ the tensor product over $R$, and
with identity object the regular $R$-bimodule ${}_RR_R$. Let
$(\cc,\Delta,\varepsilon)$ be an $R$-coring \cite{Sweedler:1975}
(i.e., a comonoid in $\bimod{R}$), we use Sweedler's notation for
the comultiplication, that is $\Delta(c)\,=\,
c_{(1)}\tensor{R}c_{(2)}$, for every $c \in \cc$ (summation is
understood). Given an $R$-bilinear morphism $\hbar:\cc\tensor{R}\cc
\to \cc\tensor{R}\cc$, we denote the image of $x\tensor{}y \in
\cc\tensor{R}\cc$ by
$\hbar(x\tensor{}y)\,=\,y^{\hbar}\tensor{}x^{\hbar}$ (summation is
understood).

Now taking into account Remark \ref{JS}, we can announce the
compatibility condition in the monoidal category $\bimod{R}$.

\begin{corollary}\label{Rings}
Let $R$ be a $\Bbbk$-algebra, and $\iota: R\to \cc$ is a rings
extension. Consider $\cc$ as an $R$-bimodule by restricting $\iota$.
Assume that this $R$-bimodule admits a structure of an $R$-coring
with comultiplication and counit, respectively, $\Delta$ and
$\varepsilon$. If $\hbar:\cc\tensor{R}\cc \to \cc\tensor{R}\cc$ is a
double distributive law (i.e., an $R$-bilinear map satisfying
\eqref{A-1}-\eqref{C-4}), then the $6$-tuple
$(\cc,\Delta,\varepsilon,\mu,1_{\cc},\hbar)$ is a bi-monoid in the
monoidal category $\bimod{R}$ (Definition \ref{bi-mono}) if and only
if
\begin{enumerate}[{R}-(a)]
\item $\Delta(1_{\cc})\,\,=\,\, 1_{\cc}\tensor{R}1_{\cc}$.
\item $\Delta(xy)\,=\,
x_{(1)}y_{(1)}{}^{\hbar}\tensor{}x_{(1)}{}^{\hbar}y_{(2)}$, for
every $x,y \in \cc$.
\item $\varepsilon(1_{\cc})\,=\, 1_R$
\item $\varepsilon(xy)\,=\, \varepsilon(x)\varepsilon(y)$, for every
$x,y \in \cc$.
\end{enumerate}
In particular $R$ is a trivial bi-monoid in $\bimod{R}$.
\end{corollary}

\subsection{Algebras and coalgebras.}

\begin{example}
Let $(B,\Delta,\varepsilon)$ be $\Bbbk$-coalgebra and $(B,\mu,1_B)$
a $\Bbbk$-algebra. Denote by $\tau:B\tensork B \to B\tensork B$ the
usual flip map i.e., $\tau(x\tensor{}y)\,=\, y\tensor{}x$, for all
$x, y \in B$. One can easily check that $\tau$ satisfies all
equalities \eqref{A-1}-\eqref{A-4} and \eqref{C-1}-\eqref{C-4} with
respect to $\Delta$, $\varepsilon$, $\mu$ and $1_B$. That is, in our
terminology, $\tau$ is a double distributive law. Therefore, by
Proposition \ref{Bimonoid}, $(B,\Delta,\varepsilon,\mu,1_B, \tau)$
is a bi-monoid in the monoidal category $\rmod{\Bbbk}$ if and only
if $B$ is a bialgebra in the usual sense \cite{Sweedler:1969}.
\end{example}

\begin{example}
Let $L$ be a $\Bbbk$-module and set $B:=\Bbbk \oplus L$. An element
belonging to $B$ will be denoted by a pair $(k,x)$ where $k \in
\Bbbk$ and $x \in L$. We consider $B$ as a $\Bbbk$-algebra with
multiplication and unit defined by
$$\mu\lr{(k,x)\tensor{}(l,y)}\,=\,(k,x) (l,y)\,\, =\,\, (kl,\, ky + ly),
\qquad 1_B=(1, 0),\quad \forall\, (k,x),\,(l,y) \in B,$$ and also as
a $\Bbbk$-coalgebra with comultiplication and counit defined by
$$ \Delta(k,x)= (k,x)\tensork (1,0) + (1,0)\tensork(0,x),\qquad
\varepsilon(k,x)=(k,0),\quad \forall\, (k,x) \in B.$$ It is well
known that $B$ is not a $\Bbbk$-bialgebra, since $\Delta$ is not a
multiplicative map. Now consider the $\Bbbk$-linear map $$
\xymatrix@R=0pt{ \hbar: B\tensork B \ar@{->}[rr] & & B\tensork B \\
\sum_i (k_i,x_i)\tensor{}(l_i,y_i) \ar@{|->}[rr] & &
\sum_i\lr{(l_i,y_i)\tensor{}(k_i,0) \,+\, (l_i,0)\tensor{}(0,x_i)
\,-(0,x_i)\tensor{}(0,y_i)}}$$ Let $(k,x), (l,y), (l',y') \in B$, We
claim that $\hbar$ is a double distributive law. First, we have
\begin{eqnarray*}
  \hbar((1,0)\tensor{}(l,y)) &=& (l,y)\tensor{}(1,0)  \\
   \hbar((k,x)\tensor{}(1,0))&=& (1,0)\tensor{}(k,0) + (1,0)\tensor{}(0,x)  \\
   &=& (1,0)\tensor{}(k,x)
\end{eqnarray*}
that is $\hbar$ satisfies \eqref{A-1} and \eqref{A-3}. On the other
hand
\begin{eqnarray*}
  (B\tensor{}\varepsilon) \circ \hbar((k,x)\tensor{}(l,y))
   &=& B\tensor{}\varepsilon\lr{(l,y)\tensor{}(k,0) + (l,0)\tensor{}(0,x) - (0,x)\tensor{}(0,y)} \\
   &=& (l,y)k \,\, =\,\,
   \varepsilon\tensor{}B\lr{(k,x)\tensor{}(l,y)},
\end{eqnarray*}
\begin{eqnarray*}
  (\varepsilon\tensor{}B) \circ \hbar((k,x)\tensor{}(l,y))
   &=& \varepsilon\tensor{}B\lr{(l,y)\tensor{}(k,0) + (l,0)\tensor{}(0,x) - (0,x)\tensor{}(0,y)} \\
   &=& l(k,0) + l(0,x) \,\, =\,\, l(k,x)\,\,=\,\,
   B\tensor{}\varepsilon\lr{(k,x)\tensor{}(l,y)}.
\end{eqnarray*}
This implies that $\hbar$ satisfies equalities \eqref{C-1} and
\eqref{C-3}. The equalities \eqref{A-2} and \eqref{A-4} are obtained
from the following two computations:
\begin{multline*}
(B\tensor{}\mu) \circ (\hbar\tensor{}B)
\circ(B\tensor{}\hbar)\lr{(k,x)\tensor{}(l',y')\tensor{}(l,y)}
\\ \,=\, (B\tensor{}\mu) \circ (\hbar\tensor{}B)\lr{(k,x)\tensor{}(l,y)\tensor{}(l',0) +
(k,x)\tensor{}(l,0)\tensor{}(0,y') -
(k,x)\tensor{}(0,y')\tensor{}(0,y) }
\\ \,=\,(B\tensor{}\mu)\left[\underset{}{} \lr{(l,y)\tensor{}(k,0) + (l,0)\tensor{}(0,x)
-(0,x)\tensor{}(0,y)}\tensor{}(l',0) + \lr{(l,0)\tensor{}(k,0) +
(l,0)\tensor{}(0,x)}\tensor{}(0,y') \right. \\ \,\,
\left.\underset{}{} - \lr{(0,y')\tensor{}(k,0)-
(0,x)\tensor{}(0,y')}\tensor{}(0,y)\right]
\\ \,=\, (l,y)\tensor{}(kl',0) + (l,0)\tensor{}(0,l'x) -(0,x)\tensor{}(0,l'y) +
(l,0)\tensor{}(0,ky') - (0,y')\tensor{}(0,ky) \\ \,=\,
(l,y)\tensor{}(kl',0) + (l,0)\tensor{}(0, l'x+ky')
-(0,l'x+ky')\tensor{}(0,y) \\ \,=\,
\hbar\lr{(kl',ky'+l'x)\tensor{}(l,y)},
\end{multline*}
\begin{multline*}
(\mu\tensor{}B) \circ(B\tensor{}\hbar) \circ
(\hbar\tensor{}B)\lr{(k,x)\tensor{}(l',y')\tensor{}(l,y)} \\ \,=\,
(\mu\tensor{}B)
\circ(B\tensor{}\hbar)\lr{(l',y')\tensor{}(k,0)\tensor{}(l,y)
+(l',0)\tensor{}(0,x)\tensor{}(l,y) -(0,x)\tensor{}(0,y')\tensor{}(l,y)} \\
\,=\,(\mu\tensor{}B)\left[\underset{}{}(l',y')\tensor{}(l,y)\tensor{}(k,0)
+ (l',0)\tensor{}(l,0)\tensor{}(0,x)
-(l',0)\tensor{}(0,x)\tensor{}(0,y)\right. \\
\,\,\left.\underset{}{}-(0,x)\tensor{}(l,0)\tensor{}(0,y')
+(0,x)\tensor{}(0,y')\tensor{}(0,y)\right] \\ \,=\,
(l'l,l'y+ly')\tensor{}(k,0) + (l'l,0)\tensor{}(0,x)
-(0,l'x)\tensor{}(0,y) -(0,xl)\tensor{}(0,y') \\ \,=\,
(l'l,l'y+ly')\tensor{}(k,0) + (l'l,0)\tensor{}(0,x)
-(0,x)\tensor{}(0,l'y+ly')\\ \,=\,
\hbar\lr{(k,x)\tensor{}(l'l,l'y+ly')}.
\end{multline*}
Equalities \eqref{C-2} and \eqref{C-4}, are obtained as follows:
\begin{multline*}
(\hbar\tensor{}B) \circ (B\tensor{}\hbar) \circ
(\Delta\tensor{}B)\lr{(k,x)\tensor{}(l,y)}  \,=\, (\hbar\tensor{}B)
\circ (B\tensor{}\hbar)\lr{ (k,x)\tensor{}(1,0)\tensor{}(l,y) +
(1,0)\tensor{}(0,x)\tensor{}(l,y)} \\ \,=\,
(\hbar\tensor{}B)\lr{(k,x)\tensor{}(l,y)\tensor{}(1,0) +
(1,0)\tensor{}(l,0)\tensor{}(0,x)-(1,0)\tensor{}(0,x)\tensor{}(0,y)}
\\ \,=\, (l,y)\tensor{}(k,0)\tensor{}(1,0) +
(l,0)\tensor{}(0,x)\tensor{}(1,0) -(0,x)\tensor{}(0,y)\tensor{}(1,0)
+ (l,0)\tensor{}(1,0)\tensor{}(0,x)
-(0,x)\tensor{}(1,0)\tensor{}(0,y)\\ \,=\, (B\tensor{}\Delta) \circ
\hbar \lr{(k,x)\tensor{}(l,y)},
\end{multline*}
\begin{multline*}
(B\tensor{}\hbar) \circ (\hbar\tensor{}B) \circ (B\tensor{}\Delta)
\lr{(k,x)\tensor{}(l,y)} \,=\, (B\tensor{}\hbar) \circ
(\hbar\tensor{}B)\lr{(k,x)\tensor{}(l,y)\tensor{}(1,0) +
(k,x)\tensor{}(1,0)\tensor{}(0,y)} \\ \,=\,
(B\tensor{}\hbar)\lr{(l,y)\tensor{}(k,0)\tensor{}(1,0) +
(l,0)\tensor{}(0,x)\tensor{}(1,0) -(0,x)\tensor{}(0,y)\tensor{}(1,0)
+ (1,0)\tensor{}(k,x)\tensor{}(0,y)} \\ \,=\,
(l,y)\tensor{}(1,0)\tensor{}(k,0) +
(l,0)\tensor{}(1,0)\tensor{}(0,x) -(0,x)\tensor{}(1,0)\tensor{}(0,y)
+ (1,0)\tensor{}(0,y)\tensor{}(k,0)
-(1,0)\tensor{}(0,x)\tensor{}(0,y) \\ \,=\, (\Delta\tensor{}B) \circ
\hbar\lr{(k,x)\tensor{}(l,y)}.
\end{multline*}
This finishes the proof of the claim. Next we show that the
$6$-tuple $(B,\Delta,\varepsilon,\mu,1_B,\hbar)$ is by Proposition
\ref{Bimonoid}, a bi-monoid in the monoidal category $\rmod{\Bbbk}$.
The equalities \ref{Bimonoid}(iii)(a), \ref{Bimonoid}(iii)(c), and
\ref{Bimonoid}(iii)(d) are easily checked, and
\ref{Bimonoid}(iii)(b) is derived from the following computation:
\begin{multline*}
(\mu\tensor{}\mu) \circ (B\tensor{}\hbar\tensor{}B) \circ
(\Delta\tensor{}\Delta)\lr{(k,x)\tensor{}(l,y)} \\ \,=\,
(\mu\tensor{}\mu) \circ
(B\tensor{}\hbar\tensor{}B)\left(\underset{}{}
(k,x)\tensor{}(1,0)\tensor{}(l,y)\tensor{}(1,0) +
(k,x)\tensor{}(1,0)\tensor{}(1,0)\tensor{}(0,y) \right. \\ \,\,
\left. \underset{}{}+
(1,0)\tensor{}(0,x)\tensor{}(l,y)\tensor{}(1,0) +
(1,0)\tensor{}(0,x)\tensor{}(1,0)\tensor{}(0,y) \right) \\ \,=\,
(\mu\tensor{}\mu)\left(\underset{}{}
(k,x)\tensor{}(l,y)\tensor{}(1,0)\tensor{}(1,0) +
(k,x)\tensor{}(1,0)\tensor{}(1,0)\tensor{}(0,y) \right. \\ \,\, +
\left. \underset{}{} (1,0)\tensor{}\lr{(l,0)\tensor{}(0,x)
-(0,x)\tensor{}(0,y)}\tensor{}(1,0) +
(1,0)\tensor{}(1,0)\tensor{}(0,x)\tensor{}(0,y) \right)
\end{multline*}
and so
\begin{eqnarray*}
 (\mu\tensor{}\mu) \circ (B\tensor{}\hbar\tensor{}B) \circ
(\Delta\tensor{}\Delta)\lr{(k,x)\tensor{}(l,y)}
   &=& (kl,ky+lx)\tensor{}(1,0) + (k,0)\tensor{}(0,y) +
(l,0)\tensor{}(0,x) \\
   &=& (kl,ky+lx)\tensor{}(1,0) + (1,0)\tensor{}(0,ky+lx) \\ &=&
   \Delta\lr{(kl,ky+lx)} \,\,=\,\, \Delta  \circ
   \mu\lr{(k,x)\tensor{}(l,y)}.
\end{eqnarray*}
\end{example}

\begin{remark}
Take $\Mm$ the monoidal category $\rmod{\Bbbk}$, let $\AA$ be a
$\Bbbk$-algebra and $\CC$ a $\Bbbk$-coalgebra. If in Proposition
\ref{product-wreath} the right $\AA$-wreath is induced by a
$\Bbbk$-algebra $\TT$, then the wreath product $\AA\tensork \TT$ is
the well known smash product $\AA\sharp_{\ttt}\TT$ as was proved in
\cite[Theorem 2.5]{Caenepeel/Ion/Militaru/Zhu:2000} (see also the
references cited there). Dually, if in Proposition
\ref{proct-cowreath} the right $\CC$-cowreath is induced by a
$\Bbbk$-coalgebra $\RR$, then the cowreath product $\CC\tensork \RR$
is the well known smash coproduct $\CC_{\rrr}\ltimes\RR$
\cite[Theorem 3.4]{Caenepeel/Ion/Militaru/Zhu:2000}. In this way,
the universal properties of smash product and smash coproduct
stated, respectively, in \cite[Proposition
2.12]{Caenepeel/Ion/Militaru/Zhu:2000} and \cite[Proposition
3.8]{Caenepeel/Ion/Militaru/Zhu:2000}, are in fact a particular
cases of Proposition \ref{Univ-Prp} and Proposition
\ref{Univ-coPrp}.

On the other hand, the factorization problem \cite[Theorem
4.5]{Caenepeel/Ion/Militaru/Zhu:2000} between two $\Bbbk$-bialgebras
can be re-formulated using double distributive law. Explicitly,
given two $\Bbbk$-bialgebras $\AA$ and $\CC$ togethers with two
$\Bbbk$-linear maps $\ccc: \AA\tensor{}\CC \to \CC\tensor{}\AA$ and
$\aaa: \CC\tensor{}\AA \to \AA\tensor{}\CC$ such that $(\CC,\ccc)
\in \rR^a_{\AA}$, $(\aaa,\CC) \in \lL^a_{\AA}$, and that $(\AA,\aaa)
\in \rR^c_{\CC}$, $(\ccc,\AA) \in \lL^c_{\CC}$. So as above
$\AA\tensor{}\CC$ is a $\Bbbk$-coalgebra and $\Bbbk$-algebra which
is not necessary a $\Bbbk$-bialgebra (this is the factorization
problem). Using Corollary \ref{Rings}, one can give as in
\cite[Theorem 4.5]{Caenepeel/Ion/Militaru/Zhu:2000} a necessary and
sufficient conditions for a double distributive law $\hbar:
(\AA\tensor{}\CC)\tensor{}(\AA\tensor{}\CC) \to
(\AA\tensor{}\CC)\tensor{}(\AA\tensor{}\CC)$ in order to get a
structure of bi-monoid on $\AA\tensor{}\CC$ in the monoidal category
$\rmod{\Bbbk}$.
\end{remark}

\providecommand{\bysame}{\leavevmode\hbox
to3em{\hrulefill}\thinspace}
\providecommand{\MR}{\relax\ifhmode\unskip\space\fi MR }
\providecommand{\MRhref}[2]{
} \providecommand{\href}[2]{#2}

\end{document}